\documentclass{proc-l}

\usepackage{amssymb}
\usepackage{bbm}

\usepackage[hidelinks]{hyperref}

\newtheorem{theorem}{Theorem}
\newtheorem{lemma}{Lemma}
\newtheorem{proposition}{Proposition}
\newtheorem{corollary}{Corollary}

\theoremstyle{remark}
\newtheorem*{remark}{Remark}

\begin{document}

\title{Explicit quadratic large sieve inequality}

\author{Zihao Liu}
\address{Department of Mathematics, University College London, London WC1H 0AY, UK}
\curraddr{Department of Mathematics, Stanford University, 380 Jane Stanford Way, Stanford, CA 94305, USA}
\email{tzliu@stanford.edu}
\thanks{}


\date{}

\dedicatory{}

%    "Communicated by" -- provide editor's name; required.
%\commby{}

\begin{abstract}
    In this article, we obtain an explicit version of Heath-Brown's large sieve inequality for quadratic characters and discuss its applications to $L$-functions and quadratic fields.
\end{abstract}

\maketitle

    \section{Introduction}

    \subsection{Motivations and applications}

    In many instances of analytic number theory (e.g. \S7, \S12, and \S17 in \cite{iwaniec_analytic_2004}), one is confronted with the task of estimating mean squares of character sums
    \begin{equation*}
        A(\mathcal S,N,a)=\sum_{\chi\in\mathcal S}\left|\sum_{n\le N}a_n\chi(n)\right|^2,
    \end{equation*}
    where $\mathcal S$ is a designated finite collection of Dirichlet characters, $N\ge1$, and $a=\{a_n\}$ is any sequence of complex numbers. By the Cauchy--Schwartz inequality, one has $A(\mathcal S,N,a)\le |\mathcal S|N\|a\|^2$, so
    \begin{equation*}
        \Delta(\mathcal S,N)=\sup_{\|a\|\ne 0}{A(\mathcal S,N,a)\over \|a\|^2}<\infty.
    \end{equation*}
    A preliminary method to improve $\Delta(\mathcal S,N)$ makes use of the orthogonality relations of Dirichlet characters. Expanding $A(\mathcal S,N,a)$ gives
    \begin{equation*}
        A(\mathcal S,N,a)
        =\sum_{n\le N}|a_n|^2\sum_\chi|\chi(n)|^2+\sum_\chi\sum_{n_1\ne n_2}a_{n_1}\overline{a_{n_2}}\chi(n_1)\overline\chi(n_2).
    \end{equation*}
    The first ``diagonal'' term is $\le|\mathcal S|\|a\|^2$, and the second ``off-diagonal'' term is usually treated with the orthogonality of characters in $\mathcal S$. When $\mathcal S$ denotes all characters modulo $q$, this method gives (see \cite[Theorem 6.2]{montgomery_topics_1971})
    \begin{equation*}
        \Delta(\mathcal S,N)\le q+N,
    \end{equation*}
    and when $\mathcal S$ denotes all primitive characters with conductor $\le Q$, we have
    \begin{equation}
        \label{eqn-largesieve}
        \Delta(\mathcal S,N)\le Q^2+N-1,
    \end{equation}
    which is known as the \textit{large sieve} inequality \cite[\S27]{davenport_multiplicative_1980}.
    
    In a fundamental paper, Heath-Brown \cite[Corollary 2]{D1995} investigated the case where $\mathcal S$ is $S(Q)$, the set of all primitive quadratic characters with conductor $\le Q$ and proved that for all $\varepsilon>0$, we have
    \begin{equation}
        \label{eqn-hb}
        A(S(Q),N,a)\ll_\varepsilon(QN)^\varepsilon(Q+N)\sum_{\substack{n_1,n_2\le N\\n_1n_2=\square}}|a_{n_1}a_{n_2}|,
    \end{equation}
    which is referred to as the \textit{quadratic large sieve} inequality. The estimate \eqref{eqn-hb} is optimal up to the magnitude of the $\ll_\varepsilon$-constant.
    
    Compared to the large sieve, the inequality \eqref{eqn-hb} achieved a power saving in $Q$-aspect by relaxing the $\ell^2$-norm $\|a\|$. When $\{a_n\}$ is supported on squarefree integers, $\sum_{n_1n_2=\square}|a_{n_1}a_{n_2}|$ coincides with $\|a\|^2$. For general $\{a_n\}$, the quantity can be conveniently estimated. For instance, according to \cite[p. 267]{D1995}, one has
    \begin{equation}
        \sum_{\substack{n_1,n_2\le N\\n_1n_2=\square}}|a_{n_1}a_{n_2}|\ll(N\log N)\max_{n\le N}|a_n|^2,
    \end{equation}
    transforming the inequality \eqref{eqn-hb} into a more applicable form.
    
    We cannot replace $\sum_{n_1n_2=\square}|a_{n_1}a_{n_2}|$ in \eqref{eqn-hb} with $\|a\|^2$ when $\{a_n\}$ is an arbitrary sequence. If $a_n=\mathbbm{1}_\square(n)$, then
    \begin{align*}
        A(S(Q),N,\mathbbm{1}_\square)
        &=\sum_{\chi\in S(Q)}\Biggl|\sum_{\substack{n^2\le N\\(n,q_\chi)=1}}1\Biggr|^2\ge\sum_{2<p\le Q}\Biggl|\sum_{\substack{n^2\le N\\n\ne p}}1\Biggr|^2\gg{QN\over\log Q},
    \end{align*}
    in which $p$ denotes a prime number. However, $\|\mathbbm{1}_\square\|^2\le\sqrt N$, so $A(S(Q),N,a)\ll_\varepsilon(QN)^\varepsilon(Q+N)\|a\|^2$ is generally false.
    
    The quadratic large sieve has produced Lindel\"of-on-average-type results for various moments of $L$-functions with quadratic twists. For instance, if $L(s,\chi)$ denotes the Dirichlet $L$-function, Heath-Brown \cite[Theorem 2]{D1995} showed that
    \begin{equation}
        \label{eqn-l4}
        \sum_{\chi\in S(Q)}|L(1/2+it,\chi)|^4\ll_\varepsilon(Q(|t|+1))^{1+\varepsilon}.
    \end{equation}
    Let $f$ be an $SL_2(\mathbb Z)$-cusp form of weight $\kappa$ that is also an eigenfunction to all Hecke operators with Fourier expansion
    \begin{equation*}
        f(z)=\sum_{n\ge1}\lambda_f(n)n^{\kappa-1\over2}e(nz),
    \end{equation*}
    where $e(s)=e^{2\pi is}$, and let $L(s,f\otimes\chi)$ be its $L$-function twisted by $\chi$:
    \begin{equation*}
        L(s,f\otimes\chi)=\sum_{n\ge1}{\lambda_f(n)\chi(n)\over n^s},
    \end{equation*}
    where $e(s)=e^{2\pi is}$. Then an argument similar to the proof of \eqref{eqn-l4} leads to
    \begin{equation}
        \label{eqn-lt2}
        \sum_{\chi\in S(Q)}|L(1/2,f\otimes\chi)|^2\ll_\varepsilon Q^{1+\varepsilon}.
    \end{equation}

    In 2010, Soundararajan and Young \cite{Soundararajan2010} proved that under GRH, when $\kappa\equiv0\pmod4$ and $S^*(Q)$ denotes the collection of $\chi(\cdot)=\left(2\over\cdot\right)\left(\cdot\over m\right)$ with $1\le m\le Q$ squarefree odd, one has an asymptotic formula:
    \begin{equation}
        \label{eqn-lt2a}
        \sum_{\chi\in S^*(Q)}L(1/2,f\otimes\chi)^2\sim C_fQ\log Q.
    \end{equation}

    The quadratic large sieve inequality \eqref{eqn-hb} has recently gained renewed interest due to the breakthrough of Li \cite{li_moments_2024}. Specifically, by improving \eqref{eqn-hb} for special choices of $a=\{a_n\}$, Li dropped the requirement of GRH for the asymptotic formula \eqref{eqn-lt2a}, which consequently turns many previously conditional results unconditional.
    
    Apart from moments of $L$-functions with quadratic twists, the quadratic large sieve has been applied to study the frequency of sign changes of Fourier coefficients of half-integral weight modular forms along fundamental discriminants. Specifically, let $k\ge2$ be an integer and $g(z)$ be a weight $k+\frac12$ Hecke eigenform for the group $\Gamma_0(4)$ with Fourier expansion
    \begin{equation*}
        g(z)=\sum_{n\ge1}c(n)n^{k-\frac12\over2}e(nz).
    \end{equation*}
    Then Lester and Radziwiłł \cite{lester_signs_2021} deduced from \eqref{eqn-lt2} that when $N_g$ is the number of sign changes of the sequence $\{c(n)\}_{n\in\operatorname{supp}(c)}$, one has
    \begin{equation}
        \label{eqn-LRsign}
        N_g\gg X^{1-\varepsilon}.
    \end{equation}
    Moreover, if one assumes GRH, this can be improved to $N_g\gg X$.

    The inequality \eqref{eqn-hb} also has applications to quadratic fields. For instance, let $f(x)\in\mathbb Z[x]$ be a fixed separable polynomial of degree $\ge2$ and $\mathcal F(N)$ be the set of reduced fractions $a/b$ such that $1\le a,b\le N$. For squarefree number $m\ge1$, define $R_f(m,N)$ to be the number of $r\in\mathcal F(N)$ such that $\mathbb Q(\sqrt{f(r)})=\mathbb Q(\sqrt m)$. Then Shparlinski \cite{shparlinski_quadratic_2024} applied Heath-Brown's square sieve \cite{heath-brown_square_1984} and the quadratic large sieve to deduce that
    \begin{equation}
        \label{eqn-qfields}
        \sum_{m\le M}\mu^2(m)R_f(m,N)\ll_\varepsilon N^{4/3+\varepsilon}M^{5/6},
    \end{equation}
    which improves the trivial bound $O(N^2)$ in the range $M\le N^{4/5}$.

    For a prime $q\equiv3\pmod4$, denote by $P_q(X)$ the number of rational primes $p\le X$ that split in $\mathbb Q(\sqrt{-q})$. Then Dunn et al. \cite{dunn_bilinear_2020} used the quadratic large sieve \eqref{eqn-hb} to show that for all $\varepsilon>0$, the asymptotic
    \begin{equation}
        \label{eqn-Dsplit}
        P_q(X)\sim\frac12\cdot{X\over\log X}
    \end{equation}
    is true for almost all $q\sim Q$ provided that $X\ge Q^\varepsilon$ and $Q\to+\infty$.

    \subsection{Statement of results}

    By the unique factorization of $n$ into $kr^2$ with $k$ squarefree and the well-known decomposition of quadratic characters \cite[\S5]{davenport_multiplicative_1980}, Heath-Brown deduced the estimate \eqref{eqn-hb} by proving the following \cite[Theorem 1]{D1995}:
    \begin{equation}
        \label{eqn-hbthm}
        \sideset{}{^*}\sum_{m\le M}\left|\sideset{}{^*}\sum_{n\le N}a_n\left(\frac nm\right)\right|^2\ll_\varepsilon(MN)^\varepsilon(M+N)\|a\|^2,
    \end{equation}
    where $\sum^*$ means a sum over squarefree odd integers.

    In this paper, we establish an explicit form of this inequality:
    \begin{theorem}
        \label{th-main}
        Let $M,N\ge1$ and $\{ a_n\}_{n\sim N}$ be a sequence of complex numbers. Define
        \begin{equation*}
            \Sigma(M,N,a)=\sideset{}{^*}\sum_{m\sim M}\left|\sideset{}{^*}\sum_{n\sim N}a_n\left(\frac nm\right)\right|^2,\quad \mathcal B(M,N)=\sup_{\|a\|\ne 0}{\Sigma(M,N,a)\over\|a\|^2}.
        \end{equation*}
        Then there exists a constant $C>0$ such that for every $\varepsilon>0$, we have
        \begin{equation}
            \label{eqn-bmn}
            \mathcal B(M,N)\le\exp_4(C\varepsilon^{-1})\cdot(MN)^\varepsilon(M+N),
        \end{equation}
        where $\exp_k$ refers to the k'th iterated exponential.
    \end{theorem}
    This improvement exhibits the explicit $\varepsilon$-dependence in the implied constants, and setting $\varepsilon=C_1/\log_4X$ ($\log_k$ refers to the k'th iterated logarithm) for some large constant $C_1>0$ allows us to replace each $X^\varepsilon$-term in the estimates \eqref{eqn-hb}, \eqref{eqn-l4}, \eqref{eqn-LRsign}, \eqref{eqn-qfields}, and \eqref{eqn-Dsplit} with $\exp(C_2\log X/\log_4X)$ for some constant $C_2>C_1$ when $X\ge\exp_4(1)$. To demonstrate this refinement, we prove the following:
    \begin{corollary}[Explicit quadratic large sieve inequality]
        \label{cor-hbexplicit}
        There exists a constant $C>0$ such that for all $Q,N\ge1$ satisfying $QN\ge\exp_4(1)$,
        \begin{equation*}
            \sum_{\chi\in S(Q)}\left|\sum_{n\le N}a_n\chi(n)\right|^2\le\exp(C(QN)/\log_4(QN))(Q+N)\sum_{\substack{n_1,n_2\le N\\n_1n_2=\square}}|a_{n_1}a_{n_2}|.
        \end{equation*}
    \end{corollary}
    \begin{corollary}[Fourth moment of Dirichlet $L$-functions]
        \label{cor-l4}
        There exists a constant $C>0$ such that when $s=\frac12+it$, $T=|t|+1$, and $Q\ge\exp_4(1)$,
        \begin{equation}
            \sum_{\chi\in S(Q)}|L(s,\chi)|^4\le \exp\{C(QT)/\log_4(QT)\}(QT).
        \end{equation}
    \end{corollary}
    \begin{remark}
        When $s=\frac12$, Shen and Stucky \cite[Corollary 1.2]{shen2024fourthmomentquadraticdirichlet} showed that the right hand side is $\asymp Q(\log Q)^{10}$.
    \end{remark}
    \begin{corollary}[Rational primes splitting in $\mathbb Q(\sqrt{-q})$]
        \label{cor-pqx}
        There exists a constant $C>0$ such that the asymptotic \eqref{eqn-Dsplit} is valid for all but $\ll Q\exp(-C\log Q/\log_4Q)$ amount of primes $q\in(Q,2Q]$ that are $\equiv3\pmod4$, provided that $Q\to+\infty$ and $X\ge\exp(3C\log Q/\log_4Q)$.
    \end{corollary}

    Although our explicit quadratic large sieve \eqref{eqn-hb} can also produce an improvement on \eqref{eqn-lt2} analogous to \autoref{cor-l4}, it is not comparable to the asymptotic obtained by Soundararajan and Young \cite{Soundararajan2010} and Xiannan Li \cite{li_moments_2024}. This is because they exploited information concerning the coefficients $\{a_n\}$. The former invoked the GRH and estimates on shifted moments of $L$-functions, and the latter utilized the multiplicativity of the coefficients and the approximate functional equations for $L$-functions. In contrast, our \autoref{th-main} and \autoref{cor-hbexplicit} make no assumptions on the coefficients $\{a_n\}$. This is particularly useful when these coefficients emerge from sieves. For instance, to obtain the estimate \eqref{eqn-qfields}, Shparlinski applied the quadratic large sieve with coefficients $\{a_n\}$ given by Legendre symbols. In \autoref{cor-pqx}, we apply \autoref{th-main} with $\{a_n\}$ supported on primes.

    Our \autoref{th-main} could serve as a baseline for further improvements. As $\exp_4(C\varepsilon^{-1})$ still grows rapidly, it would be interesting to see if one can reduce the number of iterations. If the inequality \eqref{eqn-bmn} is valid with $\exp(C\varepsilon^{-1})$ instead, then $(MN)^\varepsilon$ can be replaced with a fixed positive power of $\log(MN)$, which would be more convenient for various applications.

    \subsection{Organization of the paper}

    In \autoref{sc-outline}, \autoref{sc-lemmas}, and \autoref{sc-mainproof}, we present the outline, lemmas, and the proof for \autoref{th-main}. Then in \autoref{sc-apps}, we derive the corollaries.

    \subsection{Acknowledgements}

    The work is supported by the LMS Undergraduate Research Bursary URB-2024-41. The author would like to express his sincere gratitude towards Professor Ian Petrow for suggesting the problem and for many helpful discussions. The author also appreciates the anonymous reviewer for many useful comments.

    \subsection{Notations}

    \begin{itemize}
        \item $p$ ----- a prime number.
        \item $e(s)$ ----- the additive character $e^{2\pi is}$.
        \item $\lfloor x\rfloor$ ----- the largest $n\in\mathbb Z$ for which $n\le x<n+1$.
        \item $n\sim N$ ----- the range $N<n\le2N$.
        \item $d(n)$ ----- the divisor function $\sum_{d|n}1$.
        \item $\omega(n)$ ----- number of prime divisors of $n$. e.g. $\omega(2^2\times3)=2$.
        \item $\sum_n^*$ ----- the sum over squarefree odd integers $n$.
        \item $\|a\|=\left(\sum_n|a_n|^2\right)^{1/2}$ ----- the $\ell^2$-norm for a sequence $a=\{a_n\}$.
        \item $\exp_k,\log_k$ ----- the $k$'th iterated exponential and natural logarithm.
        \item $C,C_1,C_2,\dots$ ----- absolute constants, not necessarily the same in each occurrence.
    \end{itemize}

    Unless with a subscript, constants in $O(\cdot)$ and $\ll$ are absolute.

    \section{Outline for the proof of \autoref{th-main}}\label{sc-outline}

    As our methods are built on Heath-Brown's method of recursive estimates, we use the same notations and develop a set of lemmas similar to his. For convenience in bookkeeping, Lemma X in our paper refers to a modification of \cite[Lemma X]{D1995}, and Proposition X in this paper refers to a result not present in \cite{D1995}. Unless stated otherwise, proofs of results mentioned in this section are detailed in \autoref{sc-mainproof}.

    We study $\Sigma(M,N,a)$ via the dual sum
    \begin{equation*}
        \Sigma_1(M,N,a)=\sideset{}{^*}\sum_{m\sim M}\left|\sideset{}{^*}\sum_{n\sim N}a_n\left(\frac mn\right)\right|^2.
    \end{equation*}

    Specifically, $\Sigma$ and $\Sigma_1$ are related by
    \begin{lemma}
        \label{lm-dual}
        We have $\mathcal B(M,N)\le 2\mathcal B(N,M)$. Moreover, for every sequence $\{a_n\}_{n\sim N}$, there exists $\{a_n'\}_{n\sim N}$ such that $|a_n'|=|a_n|$ and $\Sigma(M,N,a)\le 2\Sigma_1(M,N,a')$.
    \end{lemma}

    The proof is omitted as the statement is identical to \cite[Lemma 1]{D1995}. From now on, we study $\mathcal B(M,N)$ using $\Sigma_1$.

    \subsection{Preparation for Poisson summation}
    
    Let
    \begin{equation}
        \label{eqn-wdef}
        W(x)=\exp\left(\frac83-{2\over1-(x-\frac32)^2}\right)\chi_{[\frac12,\frac52]}(x).
    \end{equation}
    Then $W(x)$ is a non-negative smooth function with compact support and $\ge1$ on $[1,2]$. Although this is different from that used by Heath-Brown \cite[p. 241]{D1995}, it still fulfills the purpose and is convenient to derive explicit estimates.
    
    For each integer $m$, let $s(m)$ denote the unique positive squarefree divisor of $m$ that makes $|m/s(m)|$ a square. Then for all $K\le M/2$, we have
    \begin{align}
        \nonumber \Sigma_1
        &\le\sideset{}{^*}\sum_{m\in\mathbb Z}W\left(\frac mM\right)\left|\sideset{}{^*}\sum_{n\sim N}a_n\left(\frac mn\right)\right|^2 \\
        \label{eqn-sm}&\le\sum_{\substack{m\text{ odd}\\s(m)>K}}W\left(\frac mM\right)\left|\sideset{}{^*}\sum_{n\sim N}a_n\left(\frac mn\right)\right|^2:=\Sigma_1^{(K)}(M,N,a).
    \end{align}
    We expand $\Sigma_1^{(K)}$ as follows:
    \begin{equation*}
        \Sigma_1^{(K)}=\sum_{\substack{m\text{ odd}\\s(m)>K}}W\left(\frac mM\right)\sideset{}{^*}\sum_{n_1,n_2\sim N}a_{n_1}\overline{a_{n_2}}\left(m\over n_1n_2\right).
    \end{equation*}

    Let $\Delta=(n_1,n_2)$. Then $\chi(m)=\left(m\over n_1n_2\right)$ is a quadratic Dirichlet character with conductor $q=n_1n_2/\Delta^2$. To utilize this relation, we partition $\Sigma_1^{(K)}$ according to the magnitude of $\Delta$. Define
    \begin{equation*}
        \Sigma_2(M,N,K,\Delta,a)=\sideset{}{^*}\sum_{\substack{n_1,n_2\sim N\\(n_1,n_2)=\Delta}}a_{n_1}\overline{a_{n_2}}\sum_{\substack{m\text{ odd}\\s(m)>K}}W\left(\frac mM\right)\left(m\over n_1n_2\right),
    \end{equation*}
    so we have
    \begin{equation*}
        \Sigma_1^{(K)}(M,N,a)=\sum_{\Delta\le2N}\Sigma_2(M,N,K,\Delta,a).
    \end{equation*}
    By making trivial estimates, we can restrict our focus to $\Sigma_2$ with small $\Delta$'s:
    \begin{lemma}
        \label{lm-s2}
        Define
        \begin{equation*}
            \mathcal B(M,N,K)=\sup_{\|a\|\ne0}{\Sigma_1^{(K)}(M,N,a)\over\|a\|^2},\quad \mathcal C(M,N,K,\Delta)=\sup_{\|a\|\ne0}{\Sigma_2(M,N,K,\Delta,a)\over\|a\|^2}.
        \end{equation*}
        Then there exists an absolute constant $C>0$ such that for all $\varepsilon>0$, $K\le M/2$, and $\Delta_0\le N$, there is some $N_1\le N/\Delta_0$ such that
        \begin{equation*}
            \mathcal B(M,N,K)\le\exp_2(C\varepsilon^{-1})\cdot N^\varepsilon\mathcal B(M,N_1,K)+\sum_{\Delta\le\Delta_0}\mathcal C(M,N,K,\Delta).
        \end{equation*}
    \end{lemma}

    To handle $\mathcal C(M,N,K,\Delta)$ with $\Delta\le\Delta_0$, we express $\Sigma_2$ as a difference:
    \begin{equation*}
        \Sigma_2(M,N,K,\Delta,a)=\Sigma_3(M,N,K,\Delta,a)-\Sigma_4(M,N,K,\Delta,a),
    \end{equation*}
    in which
    \begin{equation*}
        \Sigma_3(M,N,K,\Delta,a)=\sideset{}{^*}\sum_{\substack{n_1,n_2\sim N\\(n_1,n_2)=\Delta}}a_{n_1}\overline{a_{n_2}}\sum_{m\text{ odd}}W\left(\frac mM\right)\left(m\over n_1n_2\right)
    \end{equation*}
    and
    \begin{equation*}
        \Sigma_4(M,N,K,\Delta,a)=\sideset{}{^*}\sum_{\substack{n_1,n_2\sim N\\(n_1,n_2)=\Delta}}a_{n_1}\overline{a_{n_2}}\sum_{\substack{m\text{ odd}\\s(m)\le K}}W\left(\frac mM\right)\left(m\over n_1n_2\right).
    \end{equation*}
    \subsection{Asymptotic developments}
    Applying Poisson summation formula to $\Sigma_3$ and $\Sigma_4$ produces asymptotic formulas:
    \begin{lemma}
        \label{lm-s3}
        Let $\varepsilon>0$ be arbitrary and $M,N,K,\Delta$ satisfy
        \begin{equation*}
            N>\Delta,\quad N^2M^{-1}(MN)^\varepsilon\le K\le M(MN)^{-\varepsilon}.
        \end{equation*}
        Then $\Sigma_3=M_3+E_3\|a\|^2$, where
        \begin{equation*}
            M_3=\sideset{}{^*}\sum_{b\le K}\sqrt\frac Mb\sum_{\substack{n_1,n_2\sim N\\(n_1,n_2)=\Delta}}a_{n_1}\overline{a_{n_2}}{\varphi(q)\over2q}\kappa(\Delta,q)\left(\frac bq\right)\int_0^{+\infty}W(x^2)\mathrm dx
        \end{equation*}
        and
        \begin{equation*}
            {E_3\over\exp_2(C\varepsilon^{-1})}\le1+(MN)^{\varepsilon}\Delta\frac MNm_0\sqrt{D_1D_2}\mathcal B\left(B,{N\over D_1\Delta}\right)^{\frac12}\mathcal B\left(B,{N\over D_2\Delta}\right)^{\frac12},
        \end{equation*}
        in which
        \begin{equation*}
            q={n_1n_2\over\Delta^2}, \quad
            \kappa(\Delta,q)=\prod_{p|\Delta}\left(1-\left(\frac pq\right)p^{-\frac12}\right),\quad m_0=\min\left\{1,{N/\sqrt{MB}\over D_1D_2}\right\},
        \end{equation*}
        and $B,D_1,D_2$ are some numbers satisfying
        \begin{equation}
            \label{eqn-bdd}
            D_1,D_2\gg{1\over\log(2MN)},\quad 1\ll B\ll K.
        \end{equation}
    \end{lemma}

    \begin{remark}
        The lower bound $K\ge N^2M^{-1}(MN)^\varepsilon$ is required to make truncations during Poisson summation.
    \end{remark}

    \begin{lemma}
        \label{lm-s4}
        Let $\varepsilon>0$, $N>\Delta$, and $K\le M^{1-\varepsilon}$. Then $\Sigma_4=M_4+E_4\|a\|^2$, where
        \begin{equation*}
            M_4=\sideset{}{^*}\sum_{\substack{v\le K\\(v,2\Delta)=1}}\sqrt\frac Mv\sum_{\substack{n_1,n_2\sim N\\(n_1,n_2)=\Delta}}a_{n_1}\overline{a_{n_2}}{\varphi(q\Delta)\over2q\Delta}\left(\frac vq\right)\int_0^{+\infty}W(x^2)\mathrm dx
        \end{equation*}
        and
        \begin{equation*}
            {E_4\over\exp_2(C\varepsilon^{-1})}\le1+(MN)^\varepsilon\left(M\over BD_1D_2\right)^{\frac12}\mathcal B\left(B,{N\over D_1\Delta}\right)^{\frac12}\mathcal B\left(B,{N\over D_2\Delta}\right)^{\frac12},
        \end{equation*}
        in which $B,D_1,D_2$ are some numbers satisfying \eqref{eqn-bdd} and
        \begin{equation*}
            D_1D_2\gg\varepsilon^2 (MN)^{-\varepsilon}\Delta^{-1}M^{\frac12}B^{-\frac12}.
        \end{equation*}
    \end{lemma}
    To produce upper bounds for $\Sigma_2$ (hence for $\mathcal C(M,N,K,\Delta)$), we estimate the difference of the main terms:
    \begin{lemma}
        \label{lm-m34}
        When $1\ll \Delta\ll N$ and $0<K\le M/2$, we have
        \begin{equation*}
            M_3-M_4\le\exp_2(C\varepsilon^{-1})(MN)^\varepsilon M^{\frac12}K^{-\frac12}\mathcal B(K\Delta^2L_1,NL_2),
        \end{equation*}
        in which $L_1=C_1\varepsilon^{-1}(MN)^\varepsilon$ and $L_2=C_2\varepsilon^{-1}(MN)^\varepsilon$ for some absolute $C_1,C_2>0$.
    \end{lemma}
    \subsection{Recursive estimates}
    Combining \autoref{lm-s3}, \autoref{lm-s4}, and \autoref{lm-m34}, we find that a given estimate of $\mathcal B(M,N)$ can be used to produce bounds for $\mathcal C(M,N,K,\Delta)$. By \autoref{lm-s2}, this translates into a bound for $\mathcal B(M,N,K)$, which is also an upper bound for $\mathcal B(M,N)$ according to the inequality \eqref{eqn-sm}, so we can estimate $\mathcal B(M,N)$ recursively.

    If there is some $\xi>0$ and some decreasing $F(\cdot)\ge1$ such that
    \begin{equation}
        \label{eqn-bmn-rc}
        \mathcal B(M,N)\le F(\varepsilon)(MN)^\varepsilon(M+N^\xi)
    \end{equation}
    for all $\varepsilon>0$, then our lemmas above produce a new inequality similar to \eqref{eqn-bmn-rc} except with a smaller exponent $\xi$ on $N$:
    \begin{lemma}
        \label{lm-cmnkd}
        Under the assumption \eqref{eqn-bmn-rc}, there exist absolute $C_1,C_2>0$ such that
        \begin{align*}
            \mathcal C(M,N,K,\Delta)
            &\le\exp_2(C_1\varepsilon^{-1})F(\varepsilon)(MN)^{C_2\varepsilon} \\
            &\times \Delta^4(M+N+M^{\frac12}K^{\xi-\frac12}+M^{\frac12}NK^{-\frac12}).
        \end{align*}
    \end{lemma}

    \begin{lemma}
        \label{lm-rc0}
        Under the assumption \eqref{eqn-bmn-rc}, there exist absolute $C_1,C_2>0$ such that
        \begin{align*}
            \mathcal B(M,N)\le\exp_2(C_1\varepsilon^{-1})F(\varepsilon)(MN)^{C_2\varepsilon^{\frac12}}(M+M^{1-\xi}N^{2\xi-1}).
        \end{align*}
    \end{lemma}

    \begin{lemma}
        \label{lm-rc}
        If the assumption \eqref{eqn-bmn-rc} holds for some $\xi>1$, then there exists some absolute $P>0$ such that
        \begin{equation*}
            \mathcal B(M,N)\le\exp_2(P\varepsilon^{-2})F(\varepsilon^2 P^{-1})(MN)^\varepsilon(M+N^{2-\frac1\xi}).
        \end{equation*}
    \end{lemma}

    By expanding and applying the P\'olya-Vinogradov inequality, one has the trivial bound \cite[p. 237]{D1995}:
    \begin{equation}
        \label{eqn-bmn-trivial}
        \mathcal B(M,N)\ll M+N^2\log N,
    \end{equation}
    which means we can launch the recursion with $\xi=2$ and $F(\varepsilon)=Q\varepsilon^{-1}$ for some large $Q>0$. Therefore, it follows from \autoref{lm-rc} that when
    \begin{equation}
        \xi_r=
        \begin{cases}
            2 & r=0, \\
            2-{1\over\xi_{r-1}} & r>0
        \end{cases}
    \end{equation}
    and
    \begin{equation}
        \label{eqn-frc}
        F_r(\varepsilon)=
        \begin{cases}
            Q\varepsilon^{-1} & r=0, \\
            \exp_2(P\varepsilon^{-2})F_{r-1}(\varepsilon^2 P^{-1}) & r>0,
        \end{cases}
    \end{equation}
    the following inequality holds for all $r\in\mathbb N$:
    \begin{equation*}
        \mathcal B(M,N)\le F_r(\varepsilon)(MN)^\varepsilon(M+N^{\xi_r}).
    \end{equation*}
    By induction, one easily finds $\xi_r=1+(r+1)^{-1}$, so we have
    \begin{equation}
        \label{eqn-bmn0}
        \mathcal B(M,N)\le F_{\lfloor\varepsilon^{-1}\rfloor}(\varepsilon)(MN)^{2\varepsilon}(M+N).
    \end{equation}
    It remains to estimate $F_{\lfloor\varepsilon^{-1}\rfloor}(\varepsilon)$:
    \begin{proposition}
        \label{prop-fe}
        There exists $C_2>C_1>0$ such that for all $\varepsilon>0$ one has
        \begin{equation*}
            \exp_4(C_1\varepsilon^{-1})\le F_{\lfloor\varepsilon^{-1}\rfloor}(\varepsilon)\le\exp_4(C_2\varepsilon^{-1}).
        \end{equation*}
    \end{proposition}
    Plugging this into \eqref{eqn-bmn0} and renaming $\varepsilon$ complete the proof of \autoref{th-main}.

    \begin{remark}
        The $\exp_2(P\varepsilon^{-2})$-factor in the recursion \eqref{eqn-frc} is a consequence of estimating
        \begin{equation*}
            \sum_{n\sim N}[d(n)]^k|a_n|^2
        \end{equation*}
        by $\exp_2(C\eta^{-1})N^\eta\|a\|^2$ via the standard bound for $d(n)$. This is optimal when $a_n$ is supported on primorials on which $d(n)$ is large. From this perspective, our theorem is the best one can achieve by recursive estimates with no assumptions on $\{a_n\}$.
    \end{remark}
    \section{Auxiliary lemmas}\label{sc-lemmas}
    In this section, we state and prove results necessary to the treatments of $\Sigma_3$, $\Sigma_4$, and the recursive estimates. Comes first is a monotonicity principle for $\mathcal B(M,N)$:
    \begin{lemma}
        \label{lm-bmns}
        There exist $C,C_1,C_2>0$ such that if $M_1,M_2,N_1,N_2$ satisfy
        \begin{equation*}
            M_2\ge M_1C\log(2M_1N_1),\quad N_2\ge N_1C\log(2M_1N_1),
        \end{equation*}
        then
        \begin{equation*}
            \mathcal B(M_1,N_1)\le C_1\mathcal B(M_2,N_1),\quad\mathcal B(M_1,N_1)\le C_2\mathcal B(M_1,N_2).
        \end{equation*}
    \end{lemma}
    \begin{proof}
        See \cite[Lemma 9]{D1995}.
    \end{proof}
    Now, we develop a technical lemma used to estimate bilinear forms emerging in the course of estimating $E_3$, $E_4$, and $M_3-M_4$:
    \begin{lemma}
        \label{lm-tech}
        For complex sequences $\{a_n\}_{n\sim N},\{b_n\}_{n\sim N}$, define
        \begin{equation}
            \label{eqn-lm10}
            S(M,N,D,a,b)=\sum_{d\sim D}\sideset{}{^*}\sum_{m\sim M}\Biggl|\sideset{}{^*}\sum_{\substack{(n_1,n_2)=1\\d|n_1n_2}}a_{n_1}b_{n_2}\left(m\over n_1n_2\right)\Biggr|.
        \end{equation}
        Then we have
        \begin{equation*}
            S\le\exp_2(C\eta^{-1})(MN)^\eta(D_1D_2)^{\frac12}\mathcal B(M,N/D_1)^{\frac12}\mathcal B(M,N/D_2)^{\frac12}\cdot\|a\|\cdot\|b\|,
        \end{equation*}
        in which $D_1$ and $D_2$ satisfy
        \begin{equation*}
            {1\over\log(2MN)}\ll D_1,D_2\ll D,\quad D_1D_2\asymp{D\over\log^2(2MN)}.
        \end{equation*}
    \end{lemma}
    \begin{proof}
        As the result becomes trivial for large $D$, we assume $D\ll N$. Moreover, assume without loss of generality that $\{a_n\}$ and $\{b_n\}$ are supported on squarefree odd integers.

        Since $(n_1,n_2)=1$, $d$ factors uniquely into $d_1d_2$ such that $d_1|n_1$ and $d_2|n_2$, so
        \begin{align*}
            S
            &=\sum_{\substack{d_1,d_2\\d_1d_2\sim D}}\sideset{}{^*}\sum_{m\sim M}\Biggl|\sum_{\substack{(n_1,n_2)=1\\d_1|n_1,d_2|n_2}}a_{n_1}b_{n_2}\left(m\over n_1n_2\right)\Biggr| \\
            &=\sum_{\substack{d_1,d_2\\d_1d_2\sim D}}\sideset{}{^*}\sum_{m\sim M}\Biggl|\sum_{\substack{n_1,n_2\\d_1|n_1,d_2|n_2}}a_{n_1}b_{n_2}\left(m\over n_1n_2\right)\sum_{\delta|(n_1,n_2)}\mu(\delta)\Biggr| \\
            &\le\sum_{\delta\le2N}\sum_{\substack{d_1,d_2\\d_1d_2\sim D}}\sideset{}{^*}\sum_{m\sim M}\Biggl|\sum_{\substack{[\delta,d_1]|n_1\\ [\delta,d_2]|n_2}}a_{n_1}b_{n_2}\left(m\over n_1n_2\right)\Biggr|.
        \end{align*}
        Notice that
        \begin{align}
            \sum_{\substack{d_1,d_2\\d_1d_2\sim D}}
            \nonumber &=\sum_{\mu\le{\log D\over\log2}}\sum_{d_1\sim 2^\mu}\sum_{d_2\sim D/d_1}\le\sum_\mu\sum_{d_1\sim 2^\mu}\sum_{{D\over2^{\mu+1}}<d_2\le {D\over 2^{\mu-1}}} \\
            \label{eqn-munu}&\le\sum_\mu\sum_{d_1\sim 2^\mu}\sum_\nu\sum_{\substack{{D\over2^{\mu+1}}<d_2\le {D\over 2^{\mu-1}}\\d_2\sim 2^\nu}},
        \end{align}
        so we have
        \begin{align*}
            S\le\#\{(\mu,\nu)\}\times\sum_{\delta\le 2N}\max_{\mu,\nu}\sum_{\substack{d_1\sim 2^\mu\\d_2\sim 2^\nu}}\sideset{}{^*}\sum_{m\sim M}\Biggl|\sum_{\substack{[\delta,d_1]|n_1\\ [\delta,d_2]|n_2}}a_{n_1}b_{n_2}\left(m\over n_1n_2\right)\Biggr|.
        \end{align*}
        By the relation \eqref{eqn-munu}, we see that each $\mu$ corresponds to a uniformly bounded range of $\nu$, so $\#\{(\mu,\nu)\}=O(\log N)=O(\eta^{-1}N^{\eta/4})$, which means there exists some $P_1,P_2$ such that
        \begin{equation*}
            1\ll P_1,P_2\ll D,\quad P_1P_2\asymp D,
        \end{equation*}
        and
        \begin{equation*}
            S\ll(\log N)\sum_{\delta\le2N}\sum_{\substack{d_1\sim P_1\\d_2\sim P_2}}\sideset{}{^*}\sum_{m\sim M}\Biggl|\sum_{\substack{[\delta,d_1]|n_1\\ [\delta,d_2]|n_2}}a_{n_1}b_{n_2}\left(m\over n_1n_2\right)\Biggr|.
        \end{equation*}
        Applying the Cauchy--Schwarz inequality on the right-hand side gives
        \begin{equation}
            \label{eqn-scs}
            S\ll(\log N)\Sigma_a^{\frac12}\Sigma_b^{\frac12},
        \end{equation}
        in which
        \begin{equation*}
            \Sigma_a=\sum_{\delta\le2N}\sum_{\substack{d_1\sim P_1\\d_2\sim P_2}} \sideset{}{^*}\sum_{m\sim M}\Biggl|\sum_{\substack{n\\ [\delta,d_1]|n}}a_n\left(\frac mn\right)\Biggr|^2,
        \end{equation*}
        and $\Sigma_b$ is defined similarly. By the duality principle for bilinear forms (see \cite[p. 170]{iwaniec_analytic_2004}) and \autoref{lm-dual}, we have
        \begin{align*}
            \sum_{m\sim M}\Biggl|\sum_{\substack{n\\ [\delta,d_1]|n}}a_n\left(\frac mn\right)\Biggr|^2
            &\le\mathcal B(N/d_1,M)\sum_{\substack{n\\ [\delta,d_1]|n}}|a_n|^2 \\
            &\le 2\mathcal B(M,N/d_1)\sum_{\substack{n\\ [\delta,d_1]|n}}|a_n|^2,
        \end{align*}
        so by \autoref{lm-bmns}, there is some $C_0>0$ such that when $\theta=C_0\log(2MN)$, there is $\mathcal B(M,N/d_1)\ll\mathcal B(M,\theta N/P_1)$. Consequently, we have
        \begin{align*}
            \Sigma_a
            &\ll P_2\mathcal B(M,\theta N/P_1)\sum_{\delta\le 2N}\sum_{d_1\sim P_1}\sum_{\substack{n\\ [\delta,d_1]|n}}|a_n|^2 \\
            &=P_2\mathcal B(M,\theta N/P_1)\sum_{n\sim N}|a_n|^2\sum_{\sigma|n}\mathop{\sum_{\delta\le 2N}\sum_{d_1\sim P_1}}_{[\delta,d_1]=\sigma}1 \\
            &\le P_2\mathcal B(M,\theta N/P_1)\sum_{n\sim N}|a_n|^2\sum_{\sigma|n}3^{\omega(\sigma)} \\
            &=P_2\mathcal B(M,\theta N/P_1)\sum_{n\sim N}|a_n|^24^{\omega(n)}.
        \end{align*}
        Finally, notice that when $n\sim N$ is squarefree, $4^{\omega(n)}=[d(n)]^2\le\exp_2(C_1\eta^{-1})N^{\eta/2}$ \cite[\S 5.2]{tenenbaum_introduction_1995}, so
        \begin{equation*}
            \Sigma_a\ll\exp_2(C_1\eta^{-1}) N^{\eta/2}P_2\mathcal B(M,\theta N/P_1)\|a\|^2.
        \end{equation*}
        By the same line of reasoning, we obtain a similar bound for $\Sigma_b$. Plugging these bounds into \eqref{eqn-scs} and setting $D_1=P_1/\theta,D_2=P_2/\theta$ complete the proof.
    \end{proof}
    The next few lemmas are dedicated to Poisson summations:
    \begin{lemma}
        \label{lm-poisson}
        Let $W$ be a smooth function whose support is a compact subset of $(0,+\infty)$. Then for any positive squarefree odd $q$, we have
        \begin{equation*}
            \sum_{m\ge1}W\left(\frac mM\right)\left(\frac mq\right)={M\tau(q)\over q}\sum_{h\in\mathbb Z}\widehat W\left(h\over q/M\right)\left(\frac hq\right),
        \end{equation*}
        where $\tau(q)$ denotes the quadratic Gauss sum with modulus $q$.
    \end{lemma}
    \begin{proof}
        See \cite[Lemma 11]{D1995}.
    \end{proof}
    The following lemma is used to separate variables during the treatment of $E_3$:
    \begin{lemma}
        \label{lm-mellin}
        Let $\varrho:\mathbb R\to\mathbb R$ be a Schwartz function. Define $\varrho_+(s),\varrho_-(s)$ by
        \begin{equation*}
            \varrho_{\pm}(s)=\int_0^{+\infty}\varrho(\pm x)x^{s-1}\mathrm dx.
        \end{equation*}
        Then $\varrho_+(s),\varrho_-(s)$ are holomorphic in $\Re(s)>0$ and satisfy
        \begin{equation*}
            \varrho_+(s),\varrho_-(s)\ll_{A,\sigma} |t|^{-A}\quad (s=\sigma+it).
        \end{equation*}
        In addition, for any $\sigma>0$, we have
        \begin{equation*}
            \varrho(\pm x)={1\over2\pi i}\int_{\sigma-i\infty}^{\sigma+i\infty}x^{-s}\varrho_\pm(s)\mathrm ds,
        \end{equation*}
        and when $\sigma\to0^+$, we have
        \begin{equation*}
            \int_{-\infty}^\infty|\varrho_\pm(\sigma+it)|\mathrm dt\ll\sigma^{-1}.
        \end{equation*}
    \end{lemma}
    \begin{proof}
        Every result apart from the last one is identical to those in \cite[Lemma 12]{D1995}. For the last result, notice that integration by parts gives $\rho_{\pm}(s)\ll(\sigma+|t|^2)^{-1}$.
    \end{proof}
    The following variant of the Poisson summation formula helps identify the main terms and error terms in $\Sigma_3$ and $\Sigma_4$:
    \begin{lemma}
        \label{lm-poisson2}
        Let $\psi:\mathbb R\to\mathbb R$ be a Schwartz function and $G,H:\mathbb [0,+\infty)\to\mathbb [0,+\infty)$ be increasing functions such that
        \begin{equation*}
            |\psi(x)|\le G(A)|x|^{-A},\quad|\widehat\psi(x)|\le H(A)|x|^{-A}.
        \end{equation*}
        If $A\ge0$, $0<X_1\le X\le X_2$, $L\ge(X_2/X)^2$ and $J=\min(X/X_1,X/X_2)$, then for any integer $k\ge2$, we have
        \begin{equation}
            \label{eqn-lm13}
            \begin{aligned}
                \sum_{\substack{n\in\mathbb Z\\(n,k)=1}}\psi\left(\frac nX\right)
                &={\varphi(k)\over k}X\int_{-\infty}^{+\infty}\psi(x)\mathrm dx-\psi(0)\sum_{\substack{d|k\\d\le X_2}}\mu(d) \\
                &-\sum_{\substack{d|k\\d>X_2}}\mu(d)\frac Xd\int_{-\infty}^{+\infty}\psi(x)\mathrm dx \\
                &+\sum_{1\le|l|\le L}\sum_{\substack{d|k\\X_1<d\le X_2}}\mu(d)\frac Xd\widehat\psi\left(lX\over d\right) \\
                &+\max\{G(A+1),H(A+2)\}\cdot O(XJ^{-A}).
            \end{aligned}
        \end{equation}
    \end{lemma}
    \begin{proof}
        By tracking the implied constants in the proof of \cite[Lemma 13]{D1995}, one finds that the error term is of the form
        \begin{equation*}
            \ll G(A_1)X(X_2/X)^{1-A_1}+H(A_2)X(X/X_1)^{-A_2}+H(A_3)X(X_2/X)^{2-A_3},
        \end{equation*}
        so setting $A_1=A+1$, $A_2=A$, and $A_3=A+2$ gives the desired result.
    \end{proof}
    The next few results concern explicit estimates for smooth weight functions. Comes first is a chain rule for higher-order derivatives:

    \begin{proposition}[Fa\`a di Bruno]
        \label{prop-chain}
        Let $f(x)$ be and $g(x)$ be $n$-times differentiable. Then
        \begin{equation*}
            {\mathrm d^n\over\mathrm dx^n}(f\circ g)(x)=\sum_{\substack{r_1,r_2,\dots,r_n\ge0\\r_1+2r_2+\cdots+nr_n=n}}{n!\over r_1!r_2!\cdots r_n!}f^{(n)}(g(x))\prod_{j=1}^n\left(g^{(j)}(x)\over j!\right)^{r_j}.
        \end{equation*}
        In particular, when $g(x)$ is a quadratic polynomial, we have
        \begin{equation*}
            {\mathrm d^n\over\mathrm dx^n}(f\circ g)(x)=\sum_{0\le r\le n/2}{n!\over 2^rr!(n-2r)!}f^{(n-r)}(g(x))[g'(x)]^{n-2r}[g''(x)]^r.
        \end{equation*}
    \end{proposition}

    \begin{proof}
        For the general case, see \cite[Theorem 2]{roman_formula_1980}. The following is a proof for the case where $g(x)$ is a quadratic polynomial.

        When $n=0$, there is nothing to prove. Suppose the result holds at $n$. Since
        \begin{align*}
            {\mathrm d\over\mathrm dx}
            &\left\{(f^{(n-r)}\circ g)(x)[g'(x)]^{n-2r}[g''(x)]^r\right\} \\
            &=f^{(n+1-r)}(g(x))[g'(x)]^{n+1-2r}[g''(x)]^r \\
            &+(n-2r)f^{(n-r)}(g(x))[g'(x)]^{n-1-2r}[g''(x)]^{r+1},
        \end{align*}
        writing $n=2\mu+\nu$ for $\nu\in\{0,1\}$ gives
        \begin{align}
            \label{eqn-diff0}{\mathrm d^{n+1}\over\mathrm dx^{n+1}}(f\circ g)(x)
            &=f^{(n+1)}(g(x))[g'(x)]^{n+1}\\
            \label{eqn-diff1}&+\sum_{r=1}^\mu{n!\over 2^rr!(n-2r)!}f^{(n+1-r)}(g(x))[g'(x)]^{n+1-2r}[g''(x)]^r \\
            \label{eqn-diff2}&+\sum_{r=0}^{\mu-1}{n!\over 2^rr!(n-1-2r)!}f^{(n-r)}(g(x))[g'(x)]^{n-1-2r}[g''(x)]^{r+1} \\
            \label{eqn-diff3}&+{n!(n-2\mu)\over 2^\mu\mu!(n-2\mu)!}f^{(n-\mu)}(g(x))[g'(x)]^{n-1-2\mu}[g''(x)]^{\mu+1}.
        \end{align}
        By re-indexing the term \eqref{eqn-diff2}, we compute that the sum of \eqref{eqn-diff1} and \eqref{eqn-diff2} is
        \begin{align*}
            &\sum_{r=1}^\mu\left({n!\over 2^rr!(n-2r)!}+{n!\over2^{r-1}(r-1)!(n+1-2r)!}\right) \\
            &\times f^{(n+1-r)}(g(x))[g'(x)]^{n+1-2r}[g''(x)]^r \\
            &=\sum_{r=1}^\mu{(n+1)!\over2^rr!(n+1-2r)!}f^{(n+1-r)}(g(x))[g'(x)]^{n+1-2r}[g''(x)]^r.
        \end{align*}
        Extending this sum to $r=0$ includes the right-hand side of \eqref{eqn-diff0}. When $\nu=0$, \eqref{eqn-diff3} vanishes. When $\nu=1$, \eqref{eqn-diff3} becomes
        \begin{align*}
            &={n!\over 2^\mu\mu!}f^{(n-\mu)}(g(x))[g'(x)]^{n-1-2\mu}[g''(x)]^{\mu+1} \\
            &={n!(\mu+1)\over 2^\mu(\mu+1)!}f^{(n+1-(\mu+1))}(g(x))[g'(x)]^{n+1-2(\mu+1)}[g''(x)]^{\mu+1}.
        \end{align*}
        Notice that $\mu+1=(n+1)/2=\lfloor(n+1)/2\rfloor$, so in all cases, the combination of \eqref{eqn-diff0}-\eqref{eqn-diff3} is
        \begin{equation*}
            =\sum_{0\le r\le(n+1)/2}{(n+1)!\over 2^rr!(n+1-2r)!}f^{(n+1-r)}(g(x))[g'(x)]^{n+1-2r}[g''(x)]^r,
        \end{equation*}
        completing the proof.
    \end{proof}

    \begin{proposition}
        \label{prop-wx}
        Let $W(x)$ be as in \eqref{eqn-wdef}. Then there exist absolute constants $C_1,C_2,C_3>0$ satisfying
        \begin{equation*}
            |W^{(j)}(x)|\le C_1^jj^{3j},\quad |\widehat W^{(j)}(k)|\le C_2^jC_3^AA^{3A}|k|^{-A}
        \end{equation*}
        for all $A,j\ge1$.
    \end{proposition}
    \begin{proof}
        Let $f(y)=e^y$ and $g(x)={-1\over1-x}$. Define $h(x)=f(g(x))$ and
        \begin{equation*}
            \omega(x)=e^{-{2\over1-x^2}}\chi_{[-1,1]}(x)=h(x)\cdot h(-x)\chi_{[-1,1]}(x),
        \end{equation*}
        so $W(x)=e^{\frac83}\omega\left(x-\frac32\right)$. Since for $0\le x<1$,
        \begin{equation*}
            g^{(j)}(x)=(-1)(-2)\cdots(-j)\cdot{1\over(x-1)^{j+1}}=-j![-g(x)]^{j+1}.
        \end{equation*}
        By \autoref{prop-chain}, we deduce that for $0\le x<1$,
        \begin{align*}
            |h^{(n)}(x)|
            &\le\sum_{\substack{r_1,r_2,\dots,r_n\ge0\\r_1+2r_2+\dots+nr_n=n}}{n!\over r_1!r_2!\cdots r_n!}e^{g(x)}\prod_{j=1}^n[-g(x)]^{(j+1)r_j} \\
            &\le\sum_{\substack{r_1,r_2,\dots,r_n\ge0\\r_1+2r_2+\dots+nr_n=n}}{n!\over r_1!r_2!\cdots r_n!}e^{g(x)}\prod_{j=1}^n[-g(x)]^{2jr_j} \\
            &=[-g(x)]^{2n}e^{g(x)}\sum_{\substack{r_1,r_2,\dots,r_n\ge0\\r_1+2r_2+\dots+nr_n=n}}{n!\over r_1!r_2!\cdots r_n!}\le[-g(x)]^{2n}e^{g(x)}n!e^n.
        \end{align*}
        Finally, by the elementary inequality $u^Re^{-u}\le R^Re^{-R}$ for all $R,u>0$, we have
        \begin{equation*}
            |h^{(n)}(x)|\le(2n)^{2n}n!e^{-n}\le 4^nn^{3n}.
        \end{equation*}
        By the Leibniz rule for derivatives, this means
        \begin{align*}
            |\omega^{(j)}(x)|
            &\le\sum_{n=0}^j\binom jn[4^nn^{3n}]\cdot[4^{j-n}(j-n)^{3(j-n)}] \\
            &\le 4^j\sum_{n=0}^j\binom jn[j^{3n}]\cdot[j^{3(j-n)}]\le 8^jj^{3j}.
        \end{align*}
        Thus, we can take $C_1=8e^{\frac83}$ in the estimate for $W^{(j)}(x)$. The estimate for $\widehat W^{(j)}(k)$ is a consequence of integration by parts.
    \end{proof}

    The following propositions are used to establish explicit asymptotic estimates for the approximate functional equation for $L(s,\chi)^2$:

    \begin{proposition}
        \label{prop-stirling}
        Let $s=\sigma+it$ with $\sigma$ bounded and $|t|\ge1$ and $u=x+iy$ with $x\to+\infty$. Then
        \begin{equation*}
            {\Gamma(s+u)\over\Gamma(s)}\ll C^xx^{2x}|t|^xe^{\frac\pi2|y|}.
        \end{equation*}
        uniformly for all $y\in\mathbb R$ and $|t|\ge1$.
    \end{proposition}
    \begin{proof}
        By Stirling's approximation \cite[p. 78]{titchmarsh_theory_1986}, when $|\eta|\ll |\xi|$, we have
        \begin{equation*}
            |\Gamma(\eta+i\xi)|=|\xi|^{\eta-\frac12}e^{-\frac\pi2|\xi|}e^{O(|\eta|)}.
        \end{equation*}
        Combining with the fact that $|\Gamma(\eta+i\xi)|\le\Gamma(\eta)\ll\eta^\eta$ for $\eta>0$, we have
        \begin{align*}
            {\Gamma(s+u)\over\Gamma(s)}
            &\ll C_1^xx^x{(|t+y|+1)^{\sigma+x-\frac12}e^{-\frac\pi2|t+y|}\over|t|^{\sigma-\frac12}e^{-\frac\pi2|t|}} \\
            &=C_1^xx^x|t|^{\frac12-\sigma}(|t+y|+1)^{x+\sigma-\frac12}e^{\frac\pi2(|t|-|t+y|)}.
        \end{align*}
        If $|y|\le4|t|$, it follows from $|t|-|t+y|\le|y|$ that
        \begin{equation*}
            (|t+y|+1)^{x+\sigma-\frac12}e^{\frac\pi2(|t|-|t+y|)} \ll C_2^x|t|^{x+\sigma-\frac12}e^{\frac\pi2|y|}.
        \end{equation*}
        If $|y|>4|t|$, then it follows from $|t|-|t+y|\le 2|t|-|y|<-|y|/2$ that
        \begin{align*}
            (|t+y|+1)&^{x+\sigma-\frac12}e^{\frac\pi2(|t|-|t+y|)}
            \ll C_2^x|t|^{x+\sigma-\frac12}e^{\frac\pi2|y|} \\
            &\ll C_3^x(\pi|y|/4)^{x+\sigma-\frac12}e^{-\pi|y|/4}. \\
            &\le C_3^x\left(x+\sigma-\frac12\right)^{x+\sigma-\frac12}\ll C_4^xx^x.
        \end{align*}
        Combining the bounds from both cases gives the desired result.
    \end{proof}

    \begin{proposition}
        \label{prop-vs}
        If $0\le\sigma\le1$, $|t|\ge1$, and $a\in\{0,1\}$. Define $s=\sigma+it$ and
        \begin{equation}
            \label{eqn-vs}
            V_s(y)={1\over2\pi i}\int_{3-i\infty}^{3+i\infty}(\pi y)^{-u}e^{u^2}\left\{\Gamma\left(s+u+a\over2\right)\over\Gamma\left(s+a\over2\right)\right\}^2{\mathrm du\over u}.
        \end{equation}
        Then for all $A\ge1$, one has
        \begin{equation}
            \label{eqn-vs-large}
            V_s(y)\ll C^AA^Ae^{A^2}(y/|t|)^{-A}
        \end{equation}
        as $y\to+\infty$ and
        \begin{equation}
            \label{eqn-vs-small}
            V_s(y)=1+O\{(y/|t|)^{1/3}\}
        \end{equation}
        as $y\to0^+$. The implied constants are absolute.
    \end{proposition}
    \begin{remark}
        This is an explicit form of \cite[Proposition 5.4]{iwaniec_analytic_2004} particularized for $L(s,\chi)^2$.
    \end{remark}
    \begin{proof}
        To prove the bound \eqref{eqn-vs-large}, we move the line of integration to $\Re(u)=A$, so writing $u=A+iv$, \autoref{prop-stirling} gives
        \begin{equation*}
            \left\{\Gamma\left(s+u+a\over2\right)\over\Gamma\left(s+a\over2\right)\right\}^2\ll C_1^AA^A|t|^Ae^{\frac\pi2|y|}.
        \end{equation*}
        Notice that $|e^{u^2}|=e^{A^2-v^2}$, so \eqref{eqn-vs} becomes
        \begin{equation*}
            V_s(y)\ll C_1^AA^A(\pi y)^{-A}|t|^Ae^{A^2}\int_{-\infty}^\infty{e^{\frac\pi2|v|-v^2}\over |v|+1}\mathrm dv\ll C_2^AA^Ae^{A^2}(y/|t|)^{-A}.
        \end{equation*}
        For the asymptotic \eqref{eqn-vs-small}, observe that the integrand has a simple pole at $u=0$ with residue $1$, so moving the line of integration to $\Re(u)=-\frac13$ and proceed similarly give us the desired result.
    \end{proof}

    \section{Proof of \autoref{th-main}}\label{sc-mainproof}

    This section proves the results stated in \autoref{sc-outline}. For the sake of bookkeeping, \autoref{sc-mainproof}.X gives demonstrations of lemmas and \autoref{prop-fe} in \autoref{sc-outline}.X.

    \subsection{Preparation for Poisson summation}

    \begin{proof}[Proof of \autoref{lm-s2}]
        Most of Heath-Brown's arguments can be reused. The $\varepsilon$-factor is a consequence of plugging the divisor bound $d(n)\le\exp_2(C\varepsilon^{-1})n^\varepsilon$ into
        \begin{equation*}
            \mathcal B^*(M,K)\sum_{n\sim N}d(n)^2|a_n|^2
        \end{equation*}
        of \cite[p. 254]{D1995}.
    \end{proof}

    \subsection{Asymptotic developments}\label{sc-asymp}

    In this section, we prove \autoref{lm-s3}, \autoref{lm-s4}, and \autoref{lm-m34}.
    \begin{proof}[Proof of \autoref{lm-s3}]
        According to \cite[p. 254]{D1995}, applying \autoref{lm-poisson} to $\Sigma_3$ gives
        \begin{equation*}
            \Sigma_3=\sum_{e|2\Delta}{M\mu(e)\over e}\sum_a\sum_{c\ge1}\sideset{}{^*}\sum_{\substack{n_1,n_2\sim N/\Delta\\(n_1,n_2)=1}}a_{n_1}^{(1)}\overline a_{n_2}^{(1)}\sideset{}{^*}\sum_{b\ge1}\widehat W\left(abc^2\over en_1n_2/M\right)\left(b\over n_1n_2\right),
        \end{equation*}
        in which $\sum_a$ sums over $\pm1,\pm2$ and
        \begin{equation*}
            a_n^{(1)}=a_{n\Delta}{\tau(n)\over n}\left(eac^2\over n\right).
        \end{equation*}
        Then Heath-Brown shows that when $N^2M^{-1}(MN)^\varepsilon\le K\le M(MN)^{-\varepsilon}$, the portion of $\Sigma_3$ with $b>K$ is $\ll_\varepsilon\|a\|^2$. Now, we give an explicit estimate for the $\ll_\varepsilon$ constant.

        According to \autoref{prop-wx}, we have $\widehat W(x)\ll\exp(A^2)|x|^{-A}$, so the terms with $b>K$ satisfies
        \begin{align*}
            &\ll\exp(A^2)K(\Delta KMN^{-2})^{-A}\sum_{e|2\Delta}{M\over e}\sum_a\sum_{c\ge1}\sideset{}{^*}\sum_{\substack{n_1,n_2\sim N/\Delta\\(n_1,n_2)=1}}|a_{n_1}^{(1)}\overline a_{n_2}^{(1)}|c^{-2A} \\
            &\ll\exp(A^2)KM(\Delta KMN^{-2})^{-A}\sum_{e\le 2\Delta}{1\over e}\sideset{}{^*}\sum_{\substack{n_1,n_2\sim N/\Delta}}|a_{n_1\Delta}a_{n_2\Delta}| \\
            &\ll\exp(A^2)M^2(MN)^{-A\varepsilon}{\log2\Delta\over\Delta^A}\sum_{n_1,n_2\sim N}|a_{n_1}a_{n_2}|.
        \end{align*}
        Since $|\alpha\beta|\le(|\alpha|^2+|\beta|^2)/2$, we have
        \begin{align*}
            \sum_{n_1,n_2\sim N}|a_{n_1}a_{n_2}|\le\frac12\sum_{n_1,n_2\sim N}(|a_{n_1}|^2+|a_{n_2}|^2)=N\|a\|^2,
        \end{align*}
        so setting $A=2\varepsilon^{-1}$ means the contribution of terms with $b>K$ is majorized by
        \begin{equation}
            \label{eqn-s3bk}
            \ll\exp(A^2)M^2N(MN)^{-A\varepsilon}\|a\|^2\ll\exp(4\varepsilon^{-2})\|a\|^2,
        \end{equation}
        which is sufficient for our purpose as $\exp(4\varepsilon^{-2})\ll\exp_2(\varepsilon^{-1})$.

        For terms in $\Sigma_3$ satisfying $b\le K$, Heath-Brown first divides the range $b\le K$ into $b\sim B$ (totally $O(\log M)$ subintervals). The terms of $\Sigma_3$ within $b\sim B$ can be written as
        \begin{equation*}
            \Sigma_3^{(B)}=\sum_{e|2\Delta}{M\mu(e)\over e}\sum_a\sideset{}{^*}\sum_{b\sim B}\sideset{}{^*}\sum_{\substack{n_1,n_2\sim N/\Delta\\(n_1,n_2)=1}}a_{n_1}^{(2)}\overline a_{n_2}^{(2)}S(n_1n_2),
        \end{equation*}
        where
        \begin{equation*}
            a_n^{(2)}=a_{n\Delta}{\tau(n)\over n}\left(ea\over n\right)
        \end{equation*}
        and
        \begin{equation*}
            S(q)=\frac12\sum_{\substack{c\in\mathbb Z\\(c,q)=1}}\widehat W\left(abc^2\over eq/M\right).
        \end{equation*}
        We are now in a position to apply \autoref{lm-poisson2} with
        \begin{equation}
            \label{eqn-psiw}
            X=\sqrt{eq\over Mb},\quad\psi(x)=\widehat W(ax^2),
        \end{equation}
        \begin{equation*}
            X_1=(MN)^{-\eta/2}{\sqrt e\over\Delta}N(MB)^{-\frac12},\quad X_2=(MN)^\eta{\sqrt e\over\Delta}N(MB)^{-\frac12},
        \end{equation*}
        in which $\eta>0$ is arbitrary, allowing us to take $J=(MN)^{\eta/4}$ and $L=(MN)^{3\eta}$, so $S(q)$ becomes
        \begin{equation*}
            S(q)=S_1(q)-S_2(q)-S_3(q)+S_4(q)+R(q),
        \end{equation*}
        where the $S_k(q)$ corresponds to the $k$'th main term in the expansion \eqref{eqn-lm13} and $R(q)$ corresponds to the error term.

        Using the exact arguments in \cite[p. 255]{D1995}, collecting all the $S_1$-portion in $\Sigma_3^{(B)}$ produces the main term $M_3$. Terms arising from $S_2(q)$, $S_3(q)$, and $S_4(q)$ can be bounded using Heath-Brown's original reasoning except that each application of \cite[Lemma 10]{D1995} and \cite[Lemma 12]{D1995} is replaced with our explicit \autoref{lm-tech} and \autoref{lm-mellin}. These terms eventually produce the estimate of $E_3$ in the statement of \autoref{lm-s3}.

        For the $R$-portion of $\Sigma_3^{(B)}$, we need $G(A)$ and $H(A)$ for our choice of $\psi(x)$ in \eqref{eqn-psiw}. By \autoref{prop-wx}, we can take $G(A)=C_1^AA^{3A}$ for some large fixed $C_1>0$. As for $H(A)$, plugging $f(x)=\widehat W(x)$ and $g(x)=ax^2$ into \autoref{prop-chain}, so
        \begin{align*}
            |\psi^{(n)}(x)|
            &\le\sum_{0\le r\le n/2}{n!\over2^rr!(n-2r)!}|\widehat W^{(n-r)}(ax^2)|\cdot|2ax|^{n-2r}|2a|^r \\
            &\le\sum_{0\le r\le n/2}{n!\over2^rr!(n-2r)!}|\widehat W^{(n-r)}(ax^2)|\cdot(4|x|)^n \\
            &\le (4|x|)^n\cdot\sup_{n/2\le j\le n}|\widehat W^{(j)}(ax^2)|\cdot\sum_{0\le r\le n/2}{n!\over2^rr!} \\
            &\le e^{\frac12}n!(4|x|)^n\cdot\sup_{n/2\le j\le n}|\widehat W^{(j)}(ax^2)|.
        \end{align*}
        Now, using the bounds for $\widehat W^{(j)}$ in \autoref{prop-wx} and setting $A=n+2$, we have
        \begin{equation*}
            |\psi^{(n)}(x)|\le n^{2n}C_3^AA^{3A}|x|^{-(A-n)}\le C_4^n n^{9n}|x|^{-2}.
        \end{equation*}
        Through integration by parts, we see that $H(A)=C_5^AA^{9A}$ is admissible. Similar to how we deduce the bound \eqref{eqn-s3bk}, we set $A=C_6\varepsilon^{-1}$ in the error term of \eqref{eqn-lm13}, so the error term of $\Sigma_3$ arised from $R(q)$ is majorized by $\exp_2(\varepsilon^{-1})\|a\|^2$.
    \end{proof}

    \begin{proof}[Proof of \autoref{lm-s4}]
        Similar to the proof of \autoref{lm-s3}, we apply \autoref{prop-wx} and \autoref{prop-chain} so that when $\psi(x)=W(x^2)$, we can take $G(A),H(A)=C_1^AA^{C_2A}$ for some absolute $C_1,C_2>0$. Therefore, replacing the application of \cite[Lemma 10]{D1995} and \cite[Lemma 13]{D1995} in Heath-Brown's derivation with our explicit \autoref{lm-tech} and \autoref{lm-poisson2} gives us the desired $M_4$ and $E_4$.
    \end{proof}

    \begin{proof}[Proof of \autoref{lm-m34}]
        Most of Heath-Brown's original arguments can be re-used except we apply the explicit divisor bound $d(w)\ll\exp_2(C\eta^{-1})w^\eta$ and replace the use of \cite[Lemma 10]{D1995} with our \autoref{lm-tech}.
    \end{proof}

    \subsection{Recursive estimates}\label{sc-recursive}

    In this section, we prove \autoref{lm-cmnkd}, \autoref{lm-rc0}, and \autoref{lm-rc}, which are needed to perform recursive estimates. We also derive \autoref{prop-fe}, which analyzes the growth $F_r(\varepsilon)$ mentioned in \autoref{sc-outline}.3.

    \begin{proof}[Proof of \autoref{lm-cmnkd}]
        The $F(\varepsilon)$-factor comes from \cite[(20)]{D1995}:
        \begin{align*}
            &\mathcal B\left(B,{N\over D_1\Delta}\right)^{\frac12}\mathcal B\left(B,{N\over D_2\Delta}\right)^{\frac12} \\
            &\ll F(\varepsilon)\left(B\cdot{N\over D_1}\right)^{\varepsilon/2}\left(B\cdot{N\over D_2}\right)^{\varepsilon/2}(B^\xi+N/D_1)^{\frac12}(B^\xi+N/D_2)^{\frac12} \\
            &\ll F(\varepsilon)(BN)^\varepsilon(B^\xi+N^{\frac12}B^{\xi/2}+ND_1^{-\frac12}D_2^{-\frac12}) \\
            &\ll F(\varepsilon)(MN)^\varepsilon(B^\xi+N^{\frac12}B^{\xi/2}+ND_1^{-\frac12}D_2^{-\frac12}),
        \end{align*}
        and $\mathcal B(K\Delta^2L_1,NL_2)$ term in \autoref{lm-m34} with $L_1=C_1\varepsilon^{-1}(MN)^\varepsilon$ and $L_2=C_2\varepsilon^{-1}(MN)^\varepsilon$ (see also the fifth equation in \cite[p. 264]{D1995}):
        \begin{align*}
            \mathcal B(K\Delta^2L_1,NL_2)
            &\le F(\varepsilon)(K\Delta^2L_1NL_2)^\varepsilon\{(K\Delta^2L_1)^\xi+NL_2\} \\
            &\le\exp_2(C_3\varepsilon^{-1})F(\varepsilon)(MN)^{C_4\varepsilon}\Delta^4(M^{\frac12}K^{\xi-\frac12}+M^{\frac12}NK^{-\frac12}),
        \end{align*}
        where the $\Delta^4$ factor follows from the fact that $\xi\le2$. Plugging these two explicit estimates into Heath-Brown's arguments produces the desired explicit upper bound for $\mathcal C(M,N,K,\Delta)$ in \autoref{lm-cmnkd}.
    \end{proof}

    \begin{proof}[Proof of \autoref{lm-rc0}]
        Replacing \cite[Lemma 2]{D1995} and \cite[Lemma 6]{D1995} with our explicit \autoref{lm-s2} and \autoref{lm-cmnkd}, the arguments in \cite[p. 264]{D1995} lead to an explicit form of \cite[(24)]{D1995}). That is, whenever $1\le\Delta_0<N\le N_0$,
        \begin{equation}
            \label{eqn-kmn}
            N_0(MN_0)^\varepsilon\le M,
        \end{equation}
        and $K=N_0^2M^{-1}(MN_0)^\varepsilon$, we have
        \begin{equation*}
            \mathcal B(M,N,K)\le G\{\mathcal B(M,N_1,K)+H\},
        \end{equation*}
        in which $N_1\le N\Delta_0^{-1}$,
        \begin{equation*}
            G=\exp_2(C_1\varepsilon^{-1})(MN_0)^{C_2\varepsilon}\ge2,
        \end{equation*}
        and
        \begin{equation*}
            H=F(\varepsilon)\Delta_0^5(M+M^{1-\xi}N_0^{2\xi-1}).
        \end{equation*}
        By the nature of \autoref{lm-s2}, there is a sequence $N_r$ satisfying $N_r\le N_{r-1}\Delta_0^{-1}$ and
        \begin{equation*}
            \mathcal B(M,N_{r-1},K)\le G\{\mathcal B(M,N_r,K)+H\}.
        \end{equation*}
        In other words, we have
        \begin{equation*}
            \mathcal B(M,N_{r-1},K)+{GH\over G-1}\le G\left\{\mathcal B(M,N_r,K)+{GH\over G-1}\right\},
        \end{equation*}
        so
        \begin{equation*}
            \mathcal B(M,N_0,K)\le G^r\{\mathcal B(M,N_r,K)+GH\}\le G^{2r}\{\mathcal B(M,N_r,K)+H\}.
        \end{equation*}
        Let $\Delta_0=N_0^{1/R}$, so there is some $1\le r_0\le R$ such that $N_{r_0}\le\Delta_0<N_{r_0-1}$, allowing us to conclude that
        \begin{align*}
            \mathcal B(M,N_0,K)
            &\le\exp_2(C_1\varepsilon^{-1}+\log2r_0)(MN_0)^{2C_2r_0\varepsilon}\{\mathcal B(M,N_{r_0},K) \\
            &+N_0^{5/R}F(\varepsilon)(M+M^{1-\xi}N_0^{2\xi-1})\}.
        \end{align*}
        Plugging in the trivial bound $\mathcal B(M,N_{r_0},K)\ll MN_{r_0}\le MN_0^{1/R}$ gives
        \begin{equation*}
            \mathcal B(M,N_0,K)\le\exp_2(C_1\varepsilon^{-1}+\log2R)F(\varepsilon)(MN_0)^{2C_2R\varepsilon+\frac 5R}(M+M^{1-\xi}N_0^{2\xi-1}).
        \end{equation*}
        Let $R=\varepsilon^{-\frac12}$, $C_3>2C_1$, and $C_4>2C_2$, so we have
        \begin{equation*}
            \mathcal B(M,N_0,K)\le\exp_2(C_3\varepsilon^{-1})(MN_0)^{C_4\varepsilon^{\frac12}}F(\varepsilon)(M+M^{1-\xi}N_0^{2\xi-1}).
        \end{equation*}
        Since $\Sigma_1\le\Sigma_1^{(K)}$, this is an upper bound for $\mathcal B(M,N_0)$ as well, completing the proof in the case under the assumption \eqref{eqn-kmn}. For the case where this assumption is false, it follows directly from the estimate \eqref{eqn-bmn-rc} and $\xi\le2$ that
        \begin{align*}
            \mathcal B(M,N_0)
            &\le F(\varepsilon)(MN_0)^{\varepsilon}(M+N_0^\xi) \\
            &=F(\varepsilon)(MN_0)^{\varepsilon}(M+M^{1-\xi}N_0^{2\xi-1}(M/N_0)^{\xi-1}) \\
            &\le F(\varepsilon)(MN_0)^{2\varepsilon}(M+M^{1-\xi}N_0^{2\xi-1}),
        \end{align*}
        which is consistent with the conclusion as $2\varepsilon\le C_4\varepsilon^{\frac12}$.
    \end{proof}
    \begin{proof}[Proof of \autoref{lm-rc}]
        Let $\alpha>0$ be some parameter to be determined later. If $M\le N^\alpha$, then it follows from \autoref{lm-bmns} and \autoref{lm-rc0} that there is some $C_0>0$ satisfying
        \begin{align*}
            \mathcal B(M,N)
            &\ll\mathcal B(N^\alpha C_0\varepsilon^{-\frac12}(MN)^{\varepsilon^{\frac12}},N) \\
            &\ll\exp_2(C_1\varepsilon^{-1})F(\varepsilon)[\varepsilon^{-\frac12}N^{\alpha+1}(MN)^{\varepsilon^{\frac12}}]^{C_2\varepsilon^{\frac12}}(N^\alpha+N^{\alpha(1-\xi)+2\xi-1}) \\
            &\le\exp_2(C_3\varepsilon^{-1})F(\varepsilon)(MN)^{C_4\varepsilon^{\frac12}}(N^\alpha+N^{\alpha(1-\xi)+2\xi-1}).
        \end{align*}
        Let $\alpha=2-\xi^{-1}$ so that both exponents in the last line are equal. Suppose $M\ge N^\alpha$. Then it follows from \autoref{lm-rc0} that
        \begin{align*}
            \mathcal B(M,N)
            &\le \exp_2(C_1\varepsilon^{-1})F(\varepsilon)(MN)^{C_2\varepsilon^{\frac12}}(M+N^{\alpha(1-\xi)+2\xi-1}) \\
            &=\exp_2(C_1\varepsilon^{-1})F(\varepsilon)(MN)^{C_2\varepsilon^{\frac12}}(M+N^{2-\xi^{-1}}).
        \end{align*}
        Therefore, in any case, we have
        \begin{equation*}
            \mathcal B(M,N)\le\exp_2(C_3\varepsilon^{-1})F(\varepsilon)(MN)^{C_4\varepsilon^{\frac12}}(M+N^{2-\xi^{-1}}),
        \end{equation*}
        provided that $C_3>C_1$. Let $Q=\max(C_4^2,1)$, so
        \begin{equation*}
            \mathcal B(M,N)\le\exp_2(C_3Q^2\varepsilon^{-2})F(\varepsilon^2/Q)(MN)^\varepsilon(M+N^{2-\xi^{-1}}).
        \end{equation*}
        Setting $P=C_3Q^2>Q$ and using the monotonicity of $F(\cdot)$ produce the conclusion.
    \end{proof}
    \begin{proof}[Proof of \autoref{prop-fe}]
        Define $\varepsilon_0=\varepsilon$ and $\varepsilon_r=\varepsilon_{r-1}^2/P$. Then $\varepsilon_r=P^{1-2^r}\varepsilon^{2^r}$ and it follows from the definition \eqref{eqn-frc} of $F_r(\varepsilon)$ that
        \begin{equation}
            \label{eqn-fr0}
            F_r(\varepsilon)=Q\exp(e^{1/\varepsilon_1}+e^{1/\varepsilon_2}
            +\dots+e^{1/\varepsilon_r})\varepsilon_r^{-1},
        \end{equation}
        which means when $r=\lfloor\varepsilon^{-1}\rfloor$, we have
        \begin{align*}
            F_r(\varepsilon)
            &\le Q\exp(re^{1/\varepsilon_r})\varepsilon_r^{-1}\le Q\exp\left(\varepsilon^{-1}e^{P^{2^{\varepsilon^{-1}}-1}\varepsilon^{-2^{\varepsilon^{-1}}}}\right)P^{2^{\varepsilon^{-1}}-1}\varepsilon^{-2^{\varepsilon^{-1}}} \\
            &=Q\exp_2\left(P^{2^{\varepsilon^{-1}}-1}e^{2^{\varepsilon^{-1}}\log\varepsilon^{-1}}+\log\varepsilon^{-1}\right)P^{2^{\varepsilon^{-1}}-1}e^{2^{\varepsilon^{-1}}\log\varepsilon^{-1}} \\
            &\le Q\exp_2(\exp_2(C_1\varepsilon^{-1})+\log\varepsilon^{-1})\exp_2(C_2\varepsilon^{-1})\le\exp_4(C_3\varepsilon^{-1}).
        \end{align*}
        On the other side, the expression \eqref{eqn-fr0} can be bounded below using
        \begin{align*}
            F_r(\varepsilon)
            &\ge\exp_2(1/\varepsilon_r)\ge\exp_3\{(2^r-1)(\log P+\log\varepsilon^{-1})\} \\
            &\ge\exp_3(2^rC_4)\ge\exp_3(2^{\varepsilon^{-1}-1}C_4)\ge\exp_4(C_5\varepsilon^{-1}),
        \end{align*}
        completing the proof.
    \end{proof}
    \section{Applications}\label{sc-apps}
    In this section, we apply \autoref{th-main} to derive the explicit quadratic large sieve inequality (\autoref{cor-hbexplicit}). Subsequently, we apply this result to moments of Dirichlet $L$-functions (\autoref{cor-l4}), and rational primes splitting in quadratic fields (\autoref{cor-pqx}).

    \subsection{Explicit quadratic large sieve inequality}

    To prove \autoref{cor-hbexplicit} from \autoref{th-main}, we begin by proving an explicit form of the inequality \eqref{eqn-hbthm} by breaking the range of $m$ and $n$ into dyadic intervals so that for $M,N\ge2$
    \begin{align*}
        \sideset{}{^*}\sum_{m\le M}\left|\sideset{}{^*}\sum_{n\le N}a_n\left(\frac nm\right)\right|^2
        &\ll(\log M)\max_{1\le M'\le M}\sum_{m\sim M'}\left|\sideset{}{^*}\sum_{n\le N}a_n\left(\frac nm\right)\right|^2 \\
        &\ll(\log M)(\log N)^2\max_{\substack{1\le M'\le M\\1\le N'\le N}}\sum_{m\sim M'}\left|\sideset{}{^*}\sum_{n\sim N'}a_n\left(\frac nm\right)\right|^2 \\
        &\le(\log MN)^3\max_{\substack{1\le M'\le M\\1\le N'\le N}}\mathcal B(M',N')\sum_{n\sim N'}|a_n|^2.
    \end{align*}
    Notice that $\log MN\ll\varepsilon^{-1}(MN)^{\varepsilon/2}$, so plugging in \autoref{th-main} with $\varepsilon$ replaced with $\varepsilon/2$ gives
    \begin{equation}
        \label{eqn-hbthme}
        \sideset{}{^*}\sum_{m\le M}\left|\sideset{}{^*}\sum_{n\le N}a_n\left(\frac nm\right)\right|^2\le\exp_4(C\varepsilon^{-1})(MN)^\varepsilon(M+N)\|a\|^2.
    \end{equation}

    For the inequality \eqref{eqn-hb}, we can reuse the arguments in \cite[p. 265--267]{D1995}. To obtain explicit $\varepsilon$-dependence in $\ll$-constants, it suffices to notice that all $X^\varepsilon$-terms during the derivation are introduced only through the use of \eqref{eqn-hbthme}, the divisor bound $d(n)\le\exp_2(C\varepsilon^{-1})n^\varepsilon$, and $\log X\ll\varepsilon^{-1}X^\varepsilon$. Thus, we have
    \begin{equation}
        \label{eqn-hbe}
        \sum_{\chi\in S(Q)}\left|\sum_{n\le N}a_n\chi(n)\right|^2\le\exp_4(C\varepsilon^{-1})(QN)^\varepsilon(Q+N)\sum_{\substack{n_1,n_2\le N\\n_1n_2=\square}}|a_{n_1}a_{n_2}|,
    \end{equation}
    which becomes \autoref{cor-hbexplicit} after setting $\varepsilon=C'/\log_4(QN)$ for some large $C'>0$.

    \subsection{Moments of Dirichlet $L$-functions}

    In this section, we prove \autoref{cor-l4} by showing that when $s=\frac12+it$, $T=|t|+1$, and $S'(Q)=S(2Q)\smallsetminus S(Q)$ represents all primitive quadratic characters $\chi$ modulo $q_\chi\sim Q$, one has
    \begin{equation}
        \label{eqn-l4e}
        S(Q,s)=\sum_{\chi\in S'(Q)}|L(s,\chi)|^4\ll \exp_4(C\varepsilon^{-1})(QT)^{1+\varepsilon},
    \end{equation}
    which immediately implies \eqref{eqn-l4} by dyadic summation.

    Heath-Brown's original argument began with the formula
    \begin{equation*}
        L(s,\chi)^2=\sum_{n\ge1}{d(n)\chi(n)\over n^s}e^{-n/U}-{1\over2\pi i}\int_{\alpha-i\infty}^{\alpha+i\infty}U^{w-s}\Gamma(w-s)L(w,\chi)^2\mathrm dw
    \end{equation*}
    for some $0<\alpha<\frac12$ and used a recursive estimate process to obtain the bound
    \begin{equation}
        \label{eqn-sqt}
        S(Q,\sigma+it)\ll_{\sigma,\varepsilon}(QT)^\varepsilon(Q+(QT)^{2-2\sigma})
    \end{equation}
    valid for all fixed $\sigma>\frac12$. Explicating Heath-Brown's argument indeed produces \eqref{eqn-l4e}. However, as demonstrated in \autoref{sc-recursive}, tracking $\ll$-constants in recursive estimates is a complicated procedure. In this section, we give a simpler proof of \eqref{eqn-l4e} by appealing to the approximate functional equation, which completely avoids the need for recursion.

    For primitive quadratic $\chi$ modulo $q=q_\chi$, it follows from \cite[\S5.2]{iwaniec_analytic_2004} that
    \begin{equation*}
        L(s,\chi)^2=\sum_{n\ge1}{d(n)\chi(n)\over n^{\frac12+it}}V_{\frac12+it}\left(n\over q\right)+\varepsilon(t)\sum_{n\ge1}{d(n)\chi(n)\over n^{\frac12-it}}V_{\frac12-it}\left(n\over q\right),
    \end{equation*}
    where $V_s(y)$ is defined in \eqref{eqn-vs} and $\varepsilon(t)$ is some function satisfying $|\varepsilon(t)|=1$.

    Applying Cauchy--Schwartz inequality, we obtain
    \begin{equation*}
        S(Q,s)\le 2S_1(Q,t)+2S_1(Q,-t),
    \end{equation*}
    where
    \begin{equation*}
        S_1(Q,t)=\sum_{\chi\in S'(Q)}\left|\sum_{n\ge1}{d(n)\chi(n)\over n^{\frac12+it}}V_{\frac12+it}\left(n\over q_\chi\right)\right|^2.
    \end{equation*}
    By \autoref{prop-vs}, when $X\gg QT$ and $A\gg1$, we have
    \begin{equation*}
        \sum_{n>X}{d(n)\chi(n)\over n^{\frac12+it}}V_{\frac12+it}\left(n\over q_\chi\right)\ll e^{2A^2}X(QTX^{-1})^A.
    \end{equation*}
    Setting $0<\delta\le1$, $X=(QT)^{1+\delta}$, and $A=2\delta^{-1}$ makes this term $\ll\exp(8\delta^{-2})$, so one has
    \begin{align*}
        S_1(Q,t)
        &\ll Q\exp(16\delta^{-2})+\sum_{\chi\in S'(Q)}\left|\sum_{n\le X}{d(n)\chi(n)\over n^{\frac12+it}}V_{\frac12+it}\left(n\over q_\chi\right)\right|^2. \\
        &\ll Q\exp(16\delta^{-2})+(\log X)^2\max_{1\le N\le X}\sum_{\chi\in S'(Q)}\left|\sum_{n\sim N}{d(n)\chi(n)\over n^{\frac12+it}}V_{\frac12+it}\left(n\over q_\chi\right)\right|^2.
    \end{align*}
    By \eqref{eqn-hbe} and the fact that $V_s(y)=O(1)$, we deduce that
    \begin{align*}
        \sum_{\chi\in S'(Q)}
        &\left|\sum_{n\sim N}{d(n)\chi(n)\over n^{\frac12+it}}V_{\frac12+it}\left(n\over q_\chi\right)\right|^2 \\
        &\ll\exp_4(C_1\varepsilon^{-1})(QN)^{\varepsilon/8}(Q+N)N\max_{n\sim N}|d(n)^2n^{-1}| \\
        &\ll\exp_4(C_2\varepsilon^{-1})(QN)^{\varepsilon/4}(Q+N)
        \le\exp_4(C_2\varepsilon^{-1})(QX)^{\varepsilon/4}(Q+X) \\
        &\ll\exp_4(C_2\varepsilon^{-1})(QT)^{1+{3\varepsilon\over4}+\delta}.
    \end{align*}
    Finally, set $\delta=\varepsilon/8$, so it follows from $\log X\ll\varepsilon^{-1}(QT)^{\varepsilon/16}$ and $\exp(16\delta^{-2})\ll\exp_4(C_2\varepsilon^{-1})$ that
    \begin{equation*}
        S_1(Q,t)\ll\exp_4(C_3\varepsilon^{-1})(QT)^{1+\varepsilon}.
    \end{equation*}
    Plugging this into $S(Q,s)$ gives the desired result.

    \subsection{Rational primes splitting in $\mathbb Q(\sqrt{-q})$}

    In this section, we prove \autoref{cor-pqx}. Let $q$ denote a prime $\equiv3\pmod4$.

    By a well known property of quadratic fields \cite[Proposition 8.5]{neukirch_algebraic_1999}, a rational prime $p$ splits in $\mathbb Q(\sqrt{-q})$ if and only if
    \begin{equation*}
        \left(-q\over p\right)=1,
    \end{equation*}
    so it follows from the reciprocity law that the number of such $p\le X$ is
    \begin{equation*}
        P_q(X)=\frac12\sum_{\substack{2<p\le X\\p\ne q}}\left[1+\left(\frac pq\right)\right].
    \end{equation*}
    Now, let $m_Q$ be the number of $q\sim Q$ such that $|P_q(X)-\frac12\pi(X)|>\delta\pi(X)$. Then
    \begin{equation*}
        m_Q\delta^2\pi(X)^2\le\frac14\sum_{q\sim Q}\left|\sum_{\substack{2<p\le X\\p\ne q}}\left(\frac pq\right)\right|^2,
    \end{equation*}
    so it follows from \autoref{th-main} that
    \begin{equation*}
        m_Q\ll\delta^{-2}\exp_4\left\{C\log(XQ)\over\log_4(XQ)\right\}\left(1+\frac QX\right).
    \end{equation*}
    Finally, we set $\log X\ge 3C\log Q/\log_4Q$ so that $\log X\ge2C\log(XQ)/\log_4(XQ)$, completing the proof of \autoref{cor-pqx}.

    %    Bibliographies can be prepared with BibTeX using amsplain,
    %    amsalpha, or (for "historical" overviews) natbib style.
    \bibliographystyle{amsplain}
    %    Insert the bibliography data here.
    \bibliography{refs}
\end{document}